\documentclass[11pt]{amsart}

\usepackage{epsfig}

\theoremstyle{plain}

\newtheorem{teo}{Theorem}
\newtheorem{cor}{Corollary}
\newtheorem{lem}{Lemma}
\newtheorem{prop}{Proposition}

\theoremstyle{remark}

\newtheorem{rem}{Remark}
\newtheorem{rems}{Remarks}
\newtheorem{nota}{Notation}
\newtheorem{example}{Example}

\theoremstyle{definition}

\newtheorem{defi}{Definition}

\newcommand{\qqed}
                {
                \qed\vspace*{3mm}
                }

\newcommand{\qedennonce}
                {
                \ \vspace*{-10pt} \\
                \qed
                \vspace*{3mm}
                }

\newcommand{\dem}
                {
                {\sl {\sc Proof.}}
                }

\newcommand{\Q}
                {
                \mathbb Q
                }

\newcommand{\Z}
                {
                {\mathbb Z}
                }

\newcommand{\geomquo}
                {
                \,/\hspace*{-0.5ex}/\,
                }

\newcommand{\compo}
                {
                {\scriptstyle \circ}
                }

\begin{document}

\title{Algebraic monoids and group embeddings}

\author{Alvaro Rittatore}
\date{\today}
\thanks{Partially supported by  a CONICYT's grant and the Universidad de
la Rep\'ublica (Uruguay).} 

\begin{abstract}

We study the geometry of algebraic monoids. We prove that
the group of  invertible elements of an irreducible algebraic 
monoid is an algebraic group, open in the monoid. Moreover, if this
group is  reductive, then the monoid is affine. We then give  a
combinatorial classification of reductive monoids by means of the
theory of spherical varieties. 
\end{abstract}

\maketitle

\section{Introduction}

Let $k$ be an algebraically closed field of arbitrary characteristic.
An {\em algebraic monoid }is an algebraic variety $S$ over $k$ with
an associative product $m:S\times S\rightarrow S$ which is a morphism
of varieties, and such that  there exists an element $1\in S$ with
$m(1,s)=m(s,1)=s$ for all $s\in S$\,. We write $m(s,s')=ss'$ for all
$s,s'\in S$\,. Let $G(S)$ be the {\em unit group }of  $S$ :
\[
G(S)=\{g\in S\ :\ \exists\, g'\in S\ ,\ gg'=g'g=1\}\ .
\]
We have  the following picture:
\[
\begin{array}{ccc}
\left\{ \mbox{\small\ algebraic monoids }\right\}&\supset&
\left\{ \mbox{\small\ algebraic groups }\right\}\\
\cup& & \cup\\
\left\{ \mbox{\small\ affine alg. monoids }\right\}&\supset&
\left\{ \mbox{\small\ affine alg. groups  }\right\}
\end{array}
\]

It is known  that if $S$ is an affine algebraic monoid, then
$G(S)$ is an algebraic group, open in $S$
(\cite{kn:putcha-monoids}).   We extend  this result to  irreducible
algebraic monoids. Moreover, the 
natural action of $G(S)\times 
G(S)$ on $S$ by left and right  multiplication has only one closed
orbit : $S$ is a {\em simple embedding }of $G(S)$\,.  
Furthermore, if we require $G(S)$ to be a reductive group, then the
monoid $S$ is affine. 
The proof of this result  is based on the 
following remark:

If $S$ is a reductive monoid, then by the Bruhat decomposition  there
exist opposite Borel subgroups $B, B^-$ of $G(S)$ such that $BB^-$ is
open in $S$\,. So $S$ is a spherical $(G(S)\times G(S))$-variety.
It follows that if $S$ has a zero, then $S$ is affine. We 
reduce the general case to that one by using  techniques of
the theory of spherical varieties. 

Moreover,  using the dictionary
between simple spherical varieties and a family of combinatorial
objects, the {\em colored cones} (\cite{kn:knop}), we classify
reductive monoids in arbitrary characteristic, in a way dual to that
one of \cite{kn:vinberg} (in characteristic zero).  

I wish to thank M. Brion, L. Renner and E. Vinberg for many useful
remarks and sugestions. 

\section{Algebraic monoids and group embeddings}

We say that an algebraic monoid $S$ is {\em irreducible }(resp. {\em
  normal}) if its 
underlying variety is irreducible (resp. normal).

An element $e\in S$ is an {\em idempotent (element) }of $S$ if
$e^2=e$\,. We note $E(S)$ the set of idempotents of $S$\,. 

We say that a monoid $S$ has a {\em zero }if there exists an element $0\in
S$ such that $0s=s0=0$ for all $s\in S$\,.

An {\em ideal }of $S$ is a non-empty subset $I\subset S$ such that
$SIS\subset I$\,. 

If there exists an ideal $Y\subset S$ contained in all the  ideals of
$S$\,, we say that $S$ has {\em kernel }$Y$\,.

If the unit group of a monoid is a reductive algebraic group, we say
that the monoid is {\em reductive.} 

A morphism $\pi : S\rightarrow S'$ between two algebraic monoids is
a {\em morphism of algebraic monoids }if $\pi(ss')=\pi(s)\pi(s')$\,,
and $\pi(1_S)=1_{S'}$\,. It is clear that  $\pi(G(S))\subset
G(S')$\,.  If $S$ and $S'$ have a zero, we impose that $\pi(0_S)=0_{S'}$\,.

These conditions are not automatically satisfied:

\begin{example}
Let $M(n)$ be the monoid of $n\times n$ matrices  with coefficients in
 $k$\,. We can identify $M(n)$ with the monoid of
endomorphisms of a $k$-vector space $V$ of dimension $n$\,. Its unit
group is $\operatorname{Gl}(n)$\,, the group of invertible matrices;
it is isomorphic to the group of automorphisms of 
$V$\,. The monoid $M(n)$ is  an affine, smooth, irreducible,
reductive monoid with zero.

The morphism $M(n)\rightarrow M(n)$\,, $A\mapsto 0$\,, commutes with
the multiplication, but it is not a morphism of  algebraic
monoids. The morphism $M(n)\rightarrow
M(n)\times M(n)$\,, $A\mapsto (A,\operatorname{id})$\,, is not a
morphism of
algebraic monoids. Indeed, the image of the zero 
matrix is not the zero of $M(n)\times M(n)$\,.
\end{example}

\begin{nota}
From now on we suppose that all the algebraic monoids  are 
irreducible, unless otherwise  mentioned.
\end{nota}

\begin{defi}

Let $G$ be an algebraic group. Consider the action of $G\times G$ on
$G$ given by 
\[
(a,b)\cdot g=agb^{-1}\quad,\quad a,b,g\in G\ .
\]

This action is transitive, and the isotropy group of $1$ is
$\Delta(G)$\,, the diagonal of $G$\,. Moreover, the action satisfies
the equation: 
\begin{equation}\label{eq:product-action}
(a,1)\cdot b = ab = (1,b^{-1})\cdot a\ .
\end{equation}

We say that a $(G\times G)$-variety $X$ is an {\em embedding of $G$
}if there exists an element $x\in X$ such that the orbit $\mathcal
O_x$ is open in $X$ and isomorphic to $G\cong (G\times
G)/\Delta(G)$\,. If there exists 
only one closed orbit, we say that $X$ is a {\em simple embedding}.
\end{defi}

\begin{teo}\label{teo:group-open}
Let $S$ be an algebraic monoid. Then $G(S)$ is a connected algebraic
group, open in $S$\,, and $S$ is a simple embedding of
$G(S)$\,. Moreover, the unique closed $(G(S)\times G(S))$-orbit in $S$
is the kernel of $S$\,.
\end{teo}

\dem  Let $m:S\times S\rightarrow S$ be the
product in $S$\,, and 
$p_1,p_2:S\times S\rightarrow S$ the projections. Then
\begin{equation}\label{inversibles_constructibles} 
G = p_1(m^{-1}(1))\cap p_2(m^{-1}(1))\ ,
\end{equation}

\noindent and $G$ is constructible in $S$\,. If we prove  that $G$
is open and dense in $S$\,, it will follow in particular  that $G$ is
connected, because $S$ is irreducible.

From equation (\ref{inversibles_constructibles}), we have
\[
G \supset p_1(m^{-1}(1)^\circ)\cap p_2(m^{-1}(1)^\circ) \ ,
\]

\noindent where $m^{-1}(1)^\circ$ is an irreducible component of
$m^{-1}(1)$ that contains $(1,1)$\,.

The product $m$ is surjective, and $S$ is irreducible, so
(\cite[prop. 6.4.5]{kn:kempf}) 
\[
\dim m^{-1}(1)^\circ \geq \dim S\ .
\]

Let $\rho_i=p_i|_{m^{-1}(1)^\circ}$\,, $i=1,2$\,. In order to show
that $G$ is dense in $S$\,, it is sufficient to prove that the
morphisms $\rho_i$ are dominant: then there exists a non-empty open
subset of $S$ contained in $\rho_1(m^{-1}(1)^\circ)\cap
\rho_2(m^{-1}(1)^\circ)\subset G$\,.

 The fibre $\rho^{-1}_1(1)$ is composed of the couples
$(x,y)\in S\times S$ such that $x=1$ and $xy=1$\,, so
$\rho_1^{-1}(1)=\{(1,1)\}$\,.

Because of the upper-semicontinuity of the function $x\mapsto \dim
\rho^{-1}(\rho(x))$\,, there exists a non-empty open subset $U\subset
m^{-1}(1)^\circ$ such that $ \dim \rho^{-1}(\rho(x))=0$ for all $x\in
U$\,. Therefore, 
\[
\dim S\leq \dim m^{-1}(1)^\circ=\dim
\overline{\rho_1(m^{-1}(1)^\circ)}\leq \dim S\ .
\]

We have then equality, and $\dim
\overline{\rho_1(m^{-1}(1)^\circ)}=\dim S$\,, and thus $\rho_1$ is a
dominant morphism. By symmetry, $\rho_2$ is a dominant morphism too.

The left multiplication $\ell_g:G\times S\rightarrow S$ by an element
$g\in G$ is 
an automorphism of $S$\,, with inverse $\ell_{g^{-1}}$\,, such that
$\ell_g(G)=G$\,. It follows immediately that $G$ is open in $S$\,, and
that every point in $G$ is smooth in $S$\,.

Next, we show that $g\mapsto g^{-1}$ is a morphism.
Before doing that, we prove that if $s\in S$ has a left inverse, then
$s$ is a unit. 

Let $s'\in S$ be such that $s's=1$\,. Then $\ell_s$ has $\ell_{s'}$ as a
left inverse, so $\ell_s$ is injective. The group $G$ being open in $S$\,,
$\ell_s(G)$ is constructible in $S$ and $\dim G=\dim S$\,. Then $\ell_s(G)$
contains a non-empty  open subset $U$ of $S$\,. But $U\cap G\neq
\emptyset$\,, so $\ell_s(G)\cap G\neq \emptyset$\,. If $g\in G\cap
\ell_s(G)$\,,  there exists $g'\in G$ such that $g=sg'$\,, and thus
$s=g{g'}^{-1}\in G$\,. Therefore, we have:
\[
m^{-1}(1)=\{ (x,y)\in S\times S\ :\ x\in G\,,\ y=x^{-1}\}\ ,
\]

\noindent and $m^{-1}(1)=m^{-1}(1)^\circ$ is irreducible, so that
$p_1|_{m^{-1}(1)}=\rho_1: m^{-1}(1) \rightarrow G$   is bijective. If
$\rho_1$ is separable, then it follows from  Zariski's Main Theorem
that $\rho_1$ is an isomorphism. Analogously, $\rho_2$ is an
isomorphism, and the inverse map $g\mapsto g^{-1}$ is obtained by the
composition $\rho_2\compo\rho_1$\,, so it is an isomorphism.

So we must prove that $\rho_1$ is separable. Every point of $G$ is smooth
in $S$\,, so we can consider the 
differential of the product $m$ at the point $(1,1)\in G\times G$\,:
\[
T_{(1,1)}m:T_{(1,1)}(S\times S)\cong T_1S\times T_1S\rightarrow
  T_1S\quad,\quad (u,v)\mapsto u+v\ .
\]

Indeed,  $m|_{\{1\}\times S}(s)=s$ and $m|_{S\times\{1\}}(s)=s$
for all $s\in S$\,; the restriction of $T_{(1,1)}m$ to the
subspaces $T_1S\times\{0\}$ and $\{0\}\times T_1S$ is then the identity.

In particular, $T_{(1,1)}m$ is surjective, and the subvariety
$m^{-1}(1)\subset S\times S$ is smooth at $(1,1)$\,, with tangent
space equal to 
\[
T_{(1,1)}m^{-1}(1)=\ker T_{(1,1)}m =\{(u,v)\in T_1S\times T_1S\ :\
u+v=0\}\ .
\] 

The differential of $\rho_1$ at the point $(1,1)$ is the restriction
to $T_{(1,1)}m^{-1}(1)$ of the differential of the projection. It follows
that $T_{(1,1)}\rho_1$ is surjective, and $\rho_1$ is separable.

We now show that $S$ is a simple embedding of $G=G(S)$\,. We consider
the action of $G\times G$ over $S$ by left and right multiplication :
\[
(g_1,g_2)\cdot s=g_1sg_2^{-1}\quad , \quad g_1,g_2\in G\ ,\ s\in S\ .
\]

The $(G\times G)$-orbit passing by $1$ is $G$\,, so $S$ is an
embedding of $G$\,. If $Y$ is a closed $(G\times G) $-orbit, then
\[
SYS=\overline{G}Y\overline{G}\subset \overline{GYG}=\overline{Y}=Y\ .
\]

It follows that every closed orbit  is an ideal. If $Y$ and $Y'$ are
closed orbits, then $\emptyset \neq YY'\subset Y\cap Y'$\,, and thus $Y=Y'$\,,
because $Y$ and $Y'$ are orbits. Let $Y$ be the unique
closed orbit. If $I\subset S$ is an ideal, then $\emptyset \neq
IY\subset Y\cap I$\,, and it follows that $Y\subset I$\,, so $Y$ is the
kernel of $S$\,.
\qqed

We remark that in the proof of the preceding theorem we have proved :

\begin{cor}
Let $S$ be an algebraic monoid. If $s\in S$ has a left (or right)
inverse, then $s$ is an unit. 
\qedennonce 
\end{cor}

On the other side, we have a partial converse to theorem
\ref{teo:group-open} (wich was proved by Vinberg in characteristic
$0$ -- see \cite{kn:vinberg}): 

\begin{prop}
Let $G$ be an algebraic group and $S$ be an affine embedding of
$G$\,. Then the product on $G$ extends to a product $\widetilde{m}:S\times
S\rightarrow S$\,, in such a way that $G(S)=G$\,.
\end{prop}

\dem 
If the product $m:G\times G\rightarrow G$ extends to a morphism
$\widetilde{m} :S\times S\rightarrow S$\,, the $\widetilde{m}$ is a
product, because $\widetilde{m}|_{G\times G}=m$ is associative. For
convenience, we note $\widetilde{m}=m$\,. It is clear that the
identity $1\in G$ satisfies $1s=s1=s$
for all $s\in S$\,. Indeed, this holds for the elements of the open
subset $G\subset S$\,. Consider the right and left actions of $G$
given by :
 \[
\begin{array}{ll}
G\times S\rightarrow S \ \ \ \ & g\cdot s=(g,e)\cdot s\\ 
S\times G\rightarrow S \ \ \ \ & s\cdot g=(e,g^{-1})\cdot s
\end{array}
\]

By construction and equation (\ref{eq:product-action}), these actions
coincide on $(G\times G)\subset (S\times S)$ with the product in
$G$\,. Then there exists a morphism 
\[
m :  (G\times S )\cup (S\times G) = U\rightarrow S
\]

\noindent which extends the product of $G$\,. Therefore, it suffices
to prove that every regular function defined in $U$ can be extended to
$S\times S$\,.

If $V$ is a vector space and $W\subset V$ a subspace, an easy
calculation shows that $W\otimes W=(V\otimes W)\cap (W\otimes
V)$\,. Applying this to $k[S]\subset k[G]$\,, we get:
\[
k[S\times S]\cong k[S]\otimes k[S] = ( k[S]\otimes k[G])\cap (
k[G]\otimes k[S])\ .
\]

On the other hand, a function $f\in k[G]\otimes k[G]$ belongs to
$k[S]\otimes k[G]$ if and only if it can be extendend to a function
defined in $S\times G$\,. The same happens for $k[G]\otimes
k[S]$\,. It follows that every function defined in $U$ can be
extended to a function defined in $S\times S$\,.

Finally, we  prove that $G(S)=G$\,. The inclusion $G\subset G(S)$ is
trivial. Let $x\in S\setminus G$ be a unit. Then
$\ell_x(x^{-1})=x\cdot x^{-1}=1\in G$\,, and there exists $g\in G$
such that $x\cdot g\in G$\,, because $G$ is open in $S$\,. It follows that
$x= (x\cdot g)\cdot g^{-1}\in G$\,.
\qqed

In order to study the geometry  of  monoids it is useful to suppose
that they are normal. The following result  allows us to doing so
without loss of generality:

\begin{lem}\label{lem:normalization} {\em (\cite{kn:renner-thesis},
\cite{kn:tesis}) }
Let $S$ be an algebraic monoid, and
$\pi:\widetilde{S}\rightarrow S$ its normalization. Then there exists 
an unique associative product $\widetilde{m}:\widetilde{S}\times
\widetilde{S}\rightarrow \widetilde{S}$ such that  $\pi\compo
\widetilde{m}=m\compo(\pi,\pi)$\,.  Moreover, the following
proprieties are verified:

{\em i) } $I\subset S$ is an ideal if and only if $\pi^{-1}(I)\subset
\widetilde{S}$ is an ideal.

{\em ii) } $S$ has a zero $0_S$ if and only if
$\pi^{-1}(0_S)=\{0_{\widetilde{S}}\}$\,, where $0_{\widetilde{S}}$ is a
zero of $\widetilde{S}$\,.

{\em iii) } $\widetilde{S}$ is an algebraic monoid, with
$G(\widetilde{S})=\pi^{-1}(G(S))$\,, and
\[
\pi|_{G(\widetilde{S})}:G(\widetilde{S})\rightarrow G(S)
\]

\noindent is an
isomorphism of groups.

{\em iv) } The restriction of $\pi$ to $E(\widetilde{S})$ gives a
bijection between $E(\widetilde{S})$ and $E(S)$\,.

It follows that $\pi$ is a morphism of (algebraic)  monoids. In
particular,  $\pi$ is equivariant for the actions of
$G(\widetilde{S})\times G(\widetilde{S})$  and $G(S)\times G(S)$\,.
\end{lem}

\dem We define the multiplication 
$\widetilde{m}:\widetilde{S}\times \widetilde{S}\rightarrow
\widetilde{S}$ by means of the universal propriety of 
normalizations: $\widetilde{m}$ is the unique morphism such that the
diagram 

\begin{center}
\begin{picture}(70,74)(0,-10)
\put(0,50){\makebox(0,0){$\widetilde{S}\times \widetilde{S}$}}
\put(0,0){\makebox(0,0){$S\times S$}}
\put(50,50){\makebox(0,0){$\widetilde{S}$}}
\put(50,0){\makebox(0,0){$S$}}
\put(-18,25){\makebox(0,0){${\scriptstyle  \pi\times \pi}$}}
\put(58,25){\makebox(0,0){${\scriptstyle \pi}$}}
\put(30,-8){\makebox(0,0){${\scriptstyle m}$}}
\put(30,58){\makebox(0,0){${\scriptstyle \widetilde{m}}$}}
\put(20,50){\vector(1,0){20}}
\put(20,0){\vector(1,0){20}}
\put(0,40){\vector(0,-1){30}}
\put(50,40){\vector(0,-1){30}}
\end{picture}
\end{center}

\noindent is commutative.

By the universal property of normalizations, we can easily prove
that $\widetilde{m}$ is associative, and by construction, $\pi$ is
then a morphism of algebraic semigroups.

We prove only iv). For the others assertions, we refer the reader to 
\cite{kn:tesis}.

It is clear that $\pi(E(\widetilde{S}))\subset E(S)$\,. Let $e\in
E(S)$\,, and $e_1\in \pi^{-1}(e)$\,. The morphisms
$\ell_{e_1e_1}=\ell_{e_1}\ell_{e_1}$ and $\ell_{e_1}$ verify the
following commutative diagram:

\begin{center}
\begin{picture}(150,80)(0,-15)
\put(0,50){\makebox(0,0){$\widetilde{S}$}}
\put(120,50){\makebox(0,0){$\widetilde{S}$}}
\put(0,0){\makebox(0,0){$S$}}
\put(120,0){\makebox(0,0){$S$}}
\put(60,50){\makebox(0,0){$\widetilde{S}$}}
\put(60,0){\makebox(0,0){$S$}}
\put(-8,25){\makebox(0,0){${\scriptstyle  \pi}$}}
\put(128,25){\makebox(0,0){${\scriptstyle  \pi}$}}
\put(68,25){\makebox(0,0){${\scriptstyle \pi}$}}
\put(30,-8){\makebox(0,0){${\scriptstyle \ell_e\ell_e=\ell_e}$}}
\put(95,-8){\makebox(0,0){${\scriptstyle \ell_e}$}}
\put(30,58){\makebox(0,0){${\scriptstyle \widetilde{\ell}_{e_1e_1}}$}}
\put(95,58){\makebox(0,0){${\scriptstyle \widetilde{\ell}_{e_1}}$}}
\put(10,50){\vector(1,0){40}}
\put(110,50){\vector(-1,0){40}}
\put(10,0){\vector(1,0){40}}
\put(110,0){\vector(-1,0){40}}
\put(0,40){\vector(0,-1){30}}
\put(120,40){\vector(0,-1){30}}
\put(60,40){\vector(0,-1){30}}
\end{picture}
\end{center}

Again by the universal property of normalizations,
$\widetilde{\ell}_{e_1}=\widetilde{\ell}_{e_1e_1}$\,. In particular,
$e_1=\widetilde{\ell}_{e_1}(1)=\widetilde{\ell}_{e_1e_1}(1)=e_1e_1$\,,
so $\pi^{-1}(E(S))\subset E(\widetilde{S})$\,. 

Let $e_1,e_2\in \pi^{-1}(e)$\,. If we replace
$\widetilde{\ell}_{e_1e_1}$ by $\widetilde{\ell}_{e_1e_2}$ in the
preceding diagram, we get again a commutative diagram, and it follows that
$e_1e_2=e_1$\,. Analogously,  $e_1e_2=e_2$\,, and thus $e_1=e_2$\,.
\qqed

\section{Reductive monoids}

We restrict ourselves to the study of reductive (irreducible)
monoids.  All the reductive groups are supposed connected.

\begin{defi}
If $G$ is a reductive group, a homogeneous space $G/H$ is {\em spherical }
if there exists a Borel subgroup $B$ of $G$ such that $BH\subset G$
is open. If $X$ is a $G$-variety with an open orbit isomorphic
to $G/H$\,, then $X$ is called an {\em  embedding of
$G/H$. }We say that $X$ is {\em simple }if the action of $G$ over $X$
has only one closed orbit. If $X$ is normal, we say that $X$ is a {\em
spherical variety.}
\end{defi}

If $G$ is a reductive group, it follows from the Bruhat
decomposition that the normal $G$-embeddings are spherical varieties,
so, by Theorem \ref{teo:group-open}, normal reductive monoids are
simple spherical varieties. The aim of this section is to study the
geometry of reductive  
monoids from the point of view of the theory of spherical varieties.
In particular, we will show that reductive monoids are affine. 

\begin{rem}
In the litterature, it is imposed to $G$-embeddings to be normal. We
don't do this in order to be able to state results valid for all reductive
monoids.
\end{rem}

We begin by proving the following 

\begin{lem}
Let $S$ be a reductive monoid with zero. Then $S$ is affine.
\end{lem}
 
\dem
More generally,  if $X$ is a simple spherical variety, with closed
orbit a point, then $X$ is affine 
(\cite[thms. 7.6 and 7.7]{kn:knop}). If $S$ is not normal, it suffices
to take the normalization $\pi:\widetilde{S}\rightarrow S$\,. It
follows from Lemma \ref{lem:normalization} that $\widetilde{S}$ is a
normal monoid with zero. So $\widetilde{S}$ is affine, and if follows
that $S$ is. Indeed, the morphism of the normalization $\pi$ is a finite
surjective  morphism, so $\widetilde{S}$ is affine if and only if $S$
is (\cite[p. 63]{kn:hartshorne-amplesubvarieties}).
\qqed

In order to continue our study, we must give {\em a priori }proofs of
results concerning idempotents in reductive monoids, which are known in
the affine case.
We recall that an {\em 1 parameter subgroup }(1-PS) of a group $G$ is a
multiplicative morphism $\lambda:k^*\rightarrow G$\,. We note
$\Xi_*(G)$ the group of 1-PS of  $G$\,.

\begin{lem}\label{lem:idempotent}
Let $S$ be a reductive monoid with unit group $G$\,, and let $T\subset
G$ be a maximal torus. Then every
$(G\times G)$-orbit contains an idempotent element belonging to
$\overline{T}$\,. 
\end{lem}

\dem
First, we show that we can approach every $(G\times G)$-orbit of $S$ by
an 1-PS of $G$\,.

If $x\in S$\,, there exists a curve germ $\xi\in
G_{((t))}=\operatorname{Hom}(\operatorname{Spec}k((t)),G)$\,, where
$k((t))$ denotes the ring of Laurent series in one variable with
coefficients in $k$\,, such that $\lim_{t\rightarrow 0}\xi(t)=x$\,. By
Iwahori's Theorem (\cite{kn:mumford}), there exists an 1-PS $\lambda$
of $G$  such that $\lambda=f\xi f'$\,, where $f$ and $f'$ belong to
$G_{k[[t]]}=\operatorname{Hom}(\operatorname{Spec}k[[t]],G)$\,. Let
$g=f(0)$ and $g'= f'(0)$\,. Then,
\[
\lim_{t\rightarrow 0} \lambda(t) =\lim_{t\rightarrow 0} f\xi(t)
f'=gxg'\in \mathcal O_x \ .
\]

The 1-PS $\lambda$ has its values in a maximal torus $T'$\,. As all
 maximal tori of $G$ are conjugated, there exists $g_0\in G$ such
that $g_0\lambda(t) g_0^{-1}\in T$ for all $t\in k^*$\,. Then
$\lambda_0(t)= g_0\lambda(t) g_0^{-1}\in \Xi_*(T)$ has a
limit in $\mathcal O_x$\,.

Next, we observe that if $\lim_{t\rightarrow 0 }\lambda_0(t)=e$\,,
then
\[
e^2=\lim_{t\rightarrow 0} \lambda(t)\lim_{t\rightarrow 0}
\lambda(t)=\lim_{t\rightarrow 0} \lambda(t)^2=\lim_{t\rightarrow 0}
\lambda(t^2)=e  \ ,
\] 

\noindent so $e$ is the idempotent we are looking for.
\qqed

\begin{lem}\label{lem:orbits-monoids}
Let $S$ be a reductive monoid with  unit group $G$\,. Let $e\in
E(S)$\,, and $Y=\mathcal O_e=GeG$\,. If $Y$ is a monoid with
identity $e$\,, then $e$ is a central element of $S$\,, and $Y$ is an
algebraic group. The isotropy group of $e$ satisfies:
\[
(G\times G)_e= (G_1\times G_1)\Delta(G)\ ,
\] 

\noindent where $G_1$ is the isotropy group of $e$ for the action of
$G$ by left multiplication. Moreover, $G_1$  is a normal
subgroup of $G$\,. In particular, $(G\times G)_e$ is a reductive group.
\end{lem}

\dem
As $e$ is the identity of $Y$\,, it follows that $ge=e(ge)=(eg)e=eg$
for all $g\in G$\,. Then $(g,g)\cdot e= geg^{-1}=gg^{-1}e=e$\,, and
$\Delta(G)\subset (G\times G)_e$\,. Moreover, $G\subset S$ is dense in
$S$\,, so  $se=es$ for all $s\in S$\,. Thus, $e$ is a
central element of $S$\,.

In particular $Y=GeG=Ge$\,, and  $\varphi_0: G/G_1\rightarrow
Y$ is a bijective morphism. Observe that the multiplication in $Y$ has
an inverse, namely $(ge)^{-1}=g^{-1}e$\,.

 If $g_1\in G_1$ and $g\in
G$\,, then 
\[
(gg_1g^{-1})\cdot e= gg_1g^{-1}e=gg_1eg^{-1}=geg^{-1}=gg^{-1}e=e
\ .
\]

It follows that $G_1$ is normal in $G$\,. Thus, $\varphi_0$ is
compatible with products.

On the other hand, it is clear that $(G_1\times
G_1)\Delta(G)\subset (G\times G)_e$\,. If $(a,b)\in (G\times G)_e$\,,
then
\[
e=(a,b)\cdot e= aeb^{-1}=ab^{-1}e=(ab^{-1},1)\cdot e\ .
\]

It follows that  $ab^{-1}\in G_1$\,, and $(a,b)=(ab^{-1},1)(b,b)\in
(G_1\times G_1)\Delta(G)$\,. 
\qqed

\begin{rem}
The multiplication in $Y$ may not be induced by that one of $G\times G$
through $G\times G\rightarrow (G\times G)/(G\times G)_e\rightarrow
Y$\,. Indeed, $(G\times G)_e$ is not always normal in $G\times G$\,.
\end{rem}

The following result of Waterhouse shows some properties of unit
groups  of affine algebraic monoids.

\begin{prop} {\rm
(\cite{kn:waterhouse},\cite[thm. 3.3.6]{kn:renner-thesis})} 
Let $G$ be an algebraic group (not necessarily reductive). Then the
following proprerties are equivalent:

i) The character group of $G$ is non trivial.

ii) The rank of the radical of $G$ is not equal to zero.

iii) There exists an affine algebraic monoid wich unit group is $G$
and such that $G\neq S$\,. 
\end{prop}

If moreover $G$ is reductive, we show that the {\em connected center
}of $S$ (i.e. the 
closure in $S$ of $Z$\,, the connected center of $G$) meets the kernel
of $S$ :

\begin{prop}\label{prop:central-idempotent}
Let $S$ be a reductive monoid with unit group $G\neq S$\,; let $Z$ be
the connected center of $G$\,, and $Y$ be the kernel of $S$\,. Then
$\overline{Z}\cap Y\neq \emptyset$\,. In particular, the
connected center of $G$ is non trivial.
\end{prop}

\dem
Let $T$ be a maximal torus of $G$\,. By Lemma \ref{lem:idempotent},
there exists an 1-PS $\lambda\in \Xi_*(T)$ whose limit exists
and is an idempotent $e\in Y\cap \overline{T}$\,.

If $W=N_G(T)/C_G(T)$ is the Weyl group of $G$\,, then $W$ acts on
$T$ by conjugation. We obtain then an action of $W$ over $\Xi_*(T)$ as
follows :
\[
w\cdot (\mu(t))=w\mu(t)w^{-1}\in \Xi_*(T)\quad \forall
\mu\in \Xi_*(T)\,,\ w\in W\ .
\]

From $W\cdot T\subset T$ we obtain $W\cdot \overline{T}\subset
\overline{T}$\,. Let $\mu=\sum_{w\in W}w\cdot
\lambda\in\Xi_*(T)$\,; it is a $W$-invariant 1-PS, with limit
$e=\prod_{w\in W}wew^{-1}\in E(Y)$\,. In this way, we have constructed
a 1-PS $\mu$ of the center of $G$ (because it is
$W$-invariant) whose limit $e$ belongs to  $Y$\,. In particular,
$\mu$ is non trivial, and $e\in Y\cap \overline{Z}$\,: the
proposition is proved. 
\qqed

It is known that every affine embedding of a semi-simple group is
 trivial (\cite{kn:putcha-monoids}). This result follows
  immediately from the  proposition above:

\begin{cor}

If the unit group of a monoid $S$ is a semi-simple group $G_0$\,,
then $G_0=S$\,.
\qedennonce
\end{cor}

\begin{teo}\label{teo:reductif-affine}

Let $S$ be a reductive monoid. Then $S$ is affine.
\end{teo}

\dem First, we consider the case where $S$ is normal. Consider the kernel
$Y$ of $S$\,, $G=G(S)$\,, $T\subset G$ a maximal 
torus and $Z\subset G$ the connected center. Let  $e\in
\overline{Z}\cap Y$ as in Proposition 
\ref{prop:central-idempotent}. Then $e$ is a central element 
of $S$\,, so $es=se$ for all $s\in S$\,. An element $y\in Y$  can be
written as a product $y=geg'=gg'e=egg'$\,. In particular, $Y$ is an
algebraic monoid with identity $e\in Y$\,.

Consider the morphism $\psi :S\rightarrow Y$\,, $s\mapsto se$\,. We
have:
\[
\left\{
\begin{array}{l}
\psi((g,g')\cdot s)=\psi(gs{g'}^{-1})=gs{g'}^{-1}e=gse{g'}^{-1}=
(g,g')\cdot(se)= (g,g')\cdot \psi(s)\\
\\
\psi(ss')=ss'e=ss'ee=ses'e= \psi(s)\psi(s')
\end{array}
\right.
\]

It follows that $\psi $ is an $(G\times G)$-equivariant morphism of
monoids. The fibre $F_e=\{s\in S\ |\ se=e\}$ is then stable under
multiplication and under the action of $(G\times G)_e$\,.

The identity $1$ belongs to $F_e$ : the fibre $F_e$ is a monoid with
zero (equal to $e$\,). If $G_1=\{g\in G\ :\ ge=e\}$\,, from
Lemma \ref{lem:orbits-monoids} we deduce  that $(G\times G)_e=(G_1\times
G_1)\Delta(G)$\,. Moreover, $F_e\cap G=G_1$\,, and it follows that
$G_1$ is an open subset of $F_e$\,, so $G_1=G(F_e)$ by Theorem
\ref{teo:group-open}. Then $F_e$ is a 
reductive monoid with zero, so it is affine.

We consider the affine variety $F=(G\times G)\times F_e$\,. The
reductive group $(G\times G)_e$ acts in $F$\,:
\[
(a,b)\cdot ((g_1,g_2),s)= ((g_1a^{-1},g_2b^{-1}),asb^{-1})\ ,
\]

\noindent where $(a,b)\in (G\times G)_e$\,, $g_1,g_2\in G$\,, $s\in
F_e$\,.  The quotient $(G\times G)*_{(G\times G)_e} F_e=
F\geomquo(G\times G)_e $ is then an affine variety. If we  prove that
this quotient is isomorphic to $S$\,, it will follow that $S$ is
affine. 

We consider the morphism
\[
\gamma: F=(G\times
G)\times F_e\rightarrow S\quad,\quad \gamma((a,b),s)=(a,b)\cdot
s=asb^{-1}\ .
\]

As $\gamma$ is constant over the $(G\times G)_e$-orbits of $F$\,, then
$\gamma$ induces a morphism $\widetilde{\gamma}:(G\times G)*_{(G\times
  G)_e}F_e\rightarrow S$\,. We show that  $\widetilde{\gamma}$ is
bijective and birational, then it follows from Zariski's Main Theorem that
$\widetilde{\gamma}$ is an isomorphism.

{\em $\widetilde{\gamma}$ is injective : } Let $((a,b),s),
((a',b'),s')\in (G\times G)\times F_e$ such that 
$\gamma((a,b),s)=\gamma((a',b'),s')$\,. Then
$(a,b)\cdot s=(a',b')\cdot s'$\,, and $s'=({a'}^{-1}a,{b'}^{-1}b)\cdot
s$\,. It follows that 
\[
(a,b)\cdot e=\psi((a,b)\cdot s)=\psi((a',b')\cdot s')=
(a',b')\cdot e\ .
\]

Then  $({a'}^{-1}a,{b'}^{-1}b)\cdot e= e$\,, and thus  
$({a'}^{-1}a,{b'}^{-1}b)\in (G\times G)_e$\,. So  
\[
\overline{((a',b'),s')} = \overline{((a',b'),({a'}^{-1}a,{b'}^{-1}b)\cdot
  s)} = \overline{((a,b),s)}\ .
\]

{\em $\widetilde{\gamma}$ is surjective : }  If $s\in S$\,, then
$\psi(s)=se=ge$ for some $g\in G$\,. Then 
$((g,1),g^{-1}s)\in (G\times G)\times F_e$\,, and
$\gamma((g,1),g^{-1}s)=s$\,.

{\em $\widetilde{\gamma}$ is birational :}

The group $G\times G$ acts on $(G\times G)*_{(G\times G)_e}F_e$ by
left multiplication in the left factor, in such a way that
$\widetilde{\gamma}$ is  $(G\times G)$-equivariant. Moreover, the
orbit $\mathcal O_{((1,1),1)}$ passing by $((1,1),1)$ is open in
$(G\times G)*_{(G\times G)_e}F_e$\,, and isomorphic to the quotient
$(G\times G) *_{(G\times G)_e} G_1$\,. The restriction of
$\widetilde{\gamma}$ 
gives then  an isomorphism   from $\mathcal O_{((1,1),1)}$ to $G$\,,
with inverse $g\mapsto \overline{((g,1),1)}\in (G\times G)*_{(G\times
G)_e}G_1$\,.

For the general case, observe that if $S$  
reductive, non-normal  monoid, its normalization $\widetilde{S}$ is a
reductive normal monoid by Lemma \ref{lem:normalization}. It follows
that $\widetilde{S}$\,,  and therefore $S$, is affine. Indeed, the
normalization $\pi:\widetilde{S}\rightarrow S$ is an  finite
surjective morphism. 

\qqed

\begin{cor}
The kernel of a normal reductive monoid is a reductive group.
\end{cor}

\dem
We have seen that if $S$ is a normal reductive monoid and $Y$ its
kernel, then 
$Y$ is a group, with identity $e\in \overline{Z}\cap Y$\,. Moreover,
if we keep the notation of the proof of theorem
\ref{teo:reductif-affine}, then 
$S\cong (G\times G)*_{(G\times G)_e} F_e$\,. It follows that 
\[
Y\cong (G\times G)\cdot \overline{((1,1),e)}\cong (G\times
G)/(G\times G)_e\ .
\]

But $(G\times G)/(G\times G)_e\cong G/G_1$\,, and $Y$ is a reductive
group. 
\qqed

\begin{rems}
{\em i) } Taking irreducible components, it is easy to prove that if
$S$ is a (reducible) reductive normal monoid, then $\overline{G(S)}$
is affine, and each unit belongs to an unique irreducible component
(\cite{kn:tesis}). Moreover, any two irreducible components of
$\overline{G(S)}$ are isomorphic as varieties.

{\em ii) } It is known that any irreducible affine monoid can be
embedded as a closed submonoid of the  monoid of endomorphisms of
some vector space $V$ of finite dimension (\cite{kn:putcha-monoids}). It
follows from Theorem \ref{teo:reductif-affine} that any reductive
monoid $S$ can be seen as the closure in $\operatorname{End}(V)$ of a
closed reductive subgroup of $\operatorname{Gl}(V)$\,, which is then
isomorphic to $G(S)$\,. 

\end{rems}

We can sumarize our results in the following

\begin{teo}
The reductive monoids are exactly the affine embeddings of reductive
groups. The commutative reductive monoids are exactly the affine
embeddings of tori.
\end{teo} 

\dem
We have only to prove the last statement. It is clear that if $S$ is a
affine embedding of a torus then $S$ is commutative. Conversely, if
$S$ is a reductive commutative monoid, then its unit group is a
reductive commutative connected group, i.e. a torus.
\qqed

Let $S$ be a reductive monoid, and $G$ its unit group. If $T$
is a maximal torus of $G$\,, then the properties of $\overline{T}$ reflect
in some sense those of $S$\,. We state first a known result on
irreducible monoids whose unit group is a torus (see for example 
(\cite{kn:oda}). 

\begin{lem}\label{lem:ideals-toric}
 Let $S$ be a  monoid whose unit group is a torus
$T$\,. Then the ideals of $S$ are exactly the 
$T$-stable closed subsets of $S$\,. In particular, every ideal of $S$ is
closed. Moreover, if $e,e'\in E(S)$\, then $ee'=e'$ if and only if
$e'\in\overline{Te}$\,.  
\qedennonce
\end{lem}

This result generalizes to reductive monoids as follows :

\begin{prop}\label{prop:ideals-idempotents}
Let $S$ be a reductive monoid, with unit group $G$\,. Then the ideals
of $S$ are exactly the 
$(G\times G)$-stable closed subsets of $S$\,. In particular all
ideals are closed.

Two idempotents $e,e'\in E(S)$ are in the same $(G\times G)$-orbit if
and only if 
they are $G$-conjugate. In particular, if $T\subset G$ is a maximal torus,
then
\[
E(S)=\bigcup_{g\in G} gE(\overline{T})g^{-1}\ .
\]
\end{prop}

\dem
 We prove the first part. For the second one, see for example
\cite[cor. 6.10]{kn:putcha-monoids}.

If $I\subset S$ is a $(G\times G)$-stable closed subset, then 
\[
SIS=\overline{G}I\overline{G}\subset \overline{GIG}=\overline{I}=I\ ,
\]

\noindent that is, $I$ is an ideal of $S$\,.

Let $I\subset S$ be an ideal of $S$\,. By definition, $I$ is
$(G\times G)$-stable. Let us prove that $I$ is closed.

From lemma \ref{lem:ideals-toric}, we get that
\[
I=\cup_{i=1}^n Ge_iG\quad , \quad e_i\in I\cap \overline{T}\cap E(S)\ .
\]

Then, its closure in $S$ is equal to $\cup_{i=1}^n
\overline{Ge_iG}$\,. We need only to prove that if $e\in E(S)\cap
I\cap\overline{T}$\,, then $\overline{GeG}\subset I$\,.

Note that $\overline{GeG}=\cup_{j=1}^{m}
G\varepsilon_jG$\,, with $\varepsilon_j\in \overline{GeG}\cap
\overline{T}$ an 
idempotent element. Moreover,
$\overline{GeG}$ is a $(G\times G)$-stable closed subset of $S$\,, so
it is an ideal of $S$\,. The intersection $GeG\cap \overline{T}$ is a
non-empty $T$-stable subset of $\overline{GeG}\cap \overline{T}$\,, so
\[
GeG\cap\overline{T}=\cup_{i=1}^l Tg_ieg'_iT = \cup_{i=1}^l Tg_ieg'_i
\quad \quad g_i,g'_i\in G 
\quad , \quad g_ieg'_i\in E(S)\ ,
\]

\noindent and the closure of $GeG\cap \overline{T}$
is an ideal of $\overline{T}$\,. Let us call this ideal $J$\,. As
$\varepsilon_j\in 
\overline{GeG}\cap \overline{T}$ for all $j=1,\ldots, 
m$\,, we get $\varepsilon_je\in J\cap E(\overline{T})$\,, and there exists
$i\in\{1,\ldots l\}$ such
that $\varepsilon_je\in  \overline {Tg_ieg'_i}$\,. It follows from lemma
\ref{lem:ideals-toric} 
that
\[
\varepsilon_je=\varepsilon_jeg_ieg'_i=(\varepsilon_jeg_i)eg'_i\in
SeS\subset I\ , 
\] 

\noindent and the proposition is proved.
\qqed

\section{Classification of reductive monoids}

In \cite{kn:vinberg}, Vinberg classified all reductive
normal monoids in characteristic zero in terms of the decomposition of
their algebra of regular functions for the action of the unit group.
This classification is dual of that of affine embeddings  
of a reductive group (see \cite{kn:vust}). From the
work of Knop (\cite{kn:knop}), we know how to classify affine embeddings of
reductive group in arbitrary characteristic. In order to do that, we
begin by recalling some results from the theory of spherical
embeddings and in particular {\em reductive embeddings
}(i.e. embeddings of reductive groups), which proofs can be found in 
\cite{kn:knop}. Next, we classify (normal) reductive monoids. In view
of Lemma \ref{lem:normalization}, we suppose that all varieties
considered are normal.

\ \\

\subsection{Reductive embeddings}
\ \\

\begin{nota}
Let  $X$ be an algebraic variety. We denote by $k(X)$ the field of rational
functions over $X$\,.
\end{nota}
\begin{defi}
Let $G$ be a reductive group, $G/H$ a spherical homogenous space, and
$X$ a simple spherical embedding of $G/H$\,, with unique closed orbit
$Y$\,. We note $\mathcal D=\mathcal D(G/H)$ 
the {\em set of colours, }i.e. the set of $B$-stable irreducible
divisors of $G/H$\,. We call the {\em (set of) colours }of $X$\,, and
we denote by $\mathcal F(X)$\,, the subset of $\mathcal D$ consisting of the
irreducible divisors such that their closure contains $Y$\,. 

We denote by $\mathcal B(X)$ the set of irreducible $G$-stable divisors of
$X$\,.

Let $\Lambda_{G/H}$ be the set of weights of $k(G/H)^{(B)}$\,, we
define
\[
{\mathcal Q}(G/H)=\operatorname{Hom}_\Z(\Lambda_{G/H}, \Q)\cong
\operatorname{Hom}_\Z(\Lambda_{G/H}, \Z)\otimes_\Z \Q\ .
\]

We call the {\em rank }of $G/H$ the dimension of the $\Q$-vector
space $\mathcal Q(G/H)$\,.

Let $M$ be an algebraic variety. A {\em valuation }of the field $k(X)$
is a function $\nu:k(X)\rightarrow \mathbb Q \cup \{-\infty\}$ with
the following properties:

{\em i) }  $\nu(f_1+f_2)\geq \min\{\nu(f_1),\nu(f_2)\}$ for all
$f_1,f_2 \in k(X)$\,.

{\em ii) } $\nu(f_1f_2)=\nu(f_1)+\nu(f_2)$ for all $f_1,f_2 \in k(X)$\,

{\em iii) }  $\nu(k^*)=0$\,, $\nu(0)=-\infty$\,.

It is easy to prove that  $\nu(k^*)=0$\,, $\nu(0)=-\infty$\,. We say
that the valuation $\nu$ is {\em normalized }if $\nu(k(X)\setminus
\{0\})=\Z$\,.

If $X$ is a $G$-variety, a valuation $\nu$ of $k(X)$ is called
$G$-invariant if $\nu(g\cdot f)=\nu(f)$ for all $g\in G$\,,
$f\in k(X)$\,, where the $G$-action over $k(X)$ is given by  $(g\cdot
f)(x)=f (g^{-1}x)$ for all $x\in X$\,. 
\end{defi}

Let $G/H$ be an spherical homogeous space, $\nu$ a  valuation
of $k(G/H)$ and  $\chi\in\Lambda_{G/H}$ a weight. If  $f_\chi, f'_\chi\in
k(G/H)^{(B)}$ are two rational functions of weight $\chi$\,,
then the function $f_\chi/f'_{\chi}$ is constant over the open
$B$-orbit\,, and thus it is constant over $G/H$\,. It follows that
$\nu(f_\chi) = \nu(f'_\chi)$\,. We  set  $\rho_\nu(\chi)=\nu(f_\chi)$\,.
If $f_\chi$ and  $f_{\chi'}$ are two rational functions of weight
$\chi$ and $\chi'$ respectively, then $f_\chi f_{\chi'}$ is a rational
function of weight $\chi+\chi'$\,, and
$\rho_\nu(\chi+\chi')=\rho_\nu(\chi)+\rho_\nu(\chi')$\,. 
Thus, every valuation  $\nu$ of $k(G/H)$  defines an  element
$\rho_v\in {\mathcal   Q}(G/H)$\,.  

If ${\mathcal V}={\mathcal V}(G/H)$ is the set of
 $G$-invariant valuations of $k(G/H)$\,, the restriction of $\rho$ to
${\mathcal V}$ gives an injection $\rho:{\mathcal V}\hookrightarrow
{\mathcal Q}(G/H)$ (\cite{kn:knop}). We will identify 
${\mathcal V}$ with its image in ${\mathcal Q}(G/H)$\,. 

The set   ${\mathcal
  V}$ is a polyhedral cone (i.e. a cone generated by a finite  number of
elements of the lattice 
$\operatorname{Hom}_\Z(\Lambda, \Z )$\,), we call it the 
{\em    valuation cone.}

Let  $D\in {\mathcal B}(X)$ be a $G$-stable irreducible divisor
 of $X$\,, and $\nu_D$ the (normalized) $G$-invariant valuation
 associated to $D$\,. 
As $D$ is determined by $\nu_D$\,, the map
\[
\rho:{\mathcal  B}(X)\rightarrow
{\mathcal V}(G/H)\quad\,,\quad \quad  D\mapsto \rho(\nu_D)
\]

\noindent is injective. We will identify  ${\mathcal B}(X)$
with its image by $\rho$\,.  

On the other hand, there is a natural mapping ${\mathcal
D}(G/H)\rightarrow {\mathcal Q}(G/H)$\,: to each color $D$  
we associate  $\rho_{\nu_D}$\,, where  $\nu_D$ is the  valuation of
$k(G/H)$ associated to $D$\,. This mapping, noted  $\rho$\,, is not
necessarily  injective (\cite{kn:knop}).

\begin{defi}
Let  $G/H$ be a spherical  homogenous space. A  {\em colored cone
  }of ${\mathcal Q}(G/H)$ is a pair $({\mathcal C},{\mathcal 
  F})$\,, where ${\mathcal F}\subset {\mathcal D}(G/H)$ is a subset
such that:

{\em i) } $\rho({\mathcal F})$  does not contain $0$\,, 

{\em ii) } ${\mathcal   C}\subset {\mathcal Q}(G/H)$ is a polyhedral
cone generated by   $\rho({\mathcal   F})$\,, and a finite number of
elements of  ${\mathcal V}$\,, 

{\em iii) } $\operatorname{int}(\mathcal C) \cap {\mathcal V}\neq
\emptyset$\,, where $\operatorname{int}(\mathcal C)$ denotes 
the relative interior of $\mathcal C$\,. 

A colored cone is {\em  strictly convex} if ${\mathcal C}$ is.

A {\em colored face }of the colored cone $(\mathcal C,\mathcal F)$ is
a colored cone   $(\tau, {\mathcal F}')$\,, with $\tau$ a face of
${\mathcal C}$\,, and ${\mathcal F}' = {\mathcal F}\cap \rho^{-1}(\tau)$\,. 
\end{defi}

\begin{prop}\label{bijection_orbites_faces}{\rm (see \cite{kn:knop})}
\indent {\em i) } There  exists a bijection between the simple
embeddings of  $G/H$ 
 and the strictly convex colored cones of  ${\mathcal
  Q}(G/H)$\,. This bijection associates to each variety $X$   the
colored cone  $({\mathcal  C}(X), {\mathcal F}(X))$\,, where
${\mathcal C}(X)$ is the  cone generated by ${\mathcal B}(X)\cup
\rho({\mathcal   F}(X))$\,. We call  $({\mathcal  
  C}(X), {\mathcal F}(X))$ the 
{\em colored cone associated to $X$\,. }We note $\mathcal
V(X)=\mathcal C(X)\cap \mathcal V(G/H)$\,.

{\em ii) } The $G$-orbits of a simple spherical variety $X$ are in
bijection with the  {\em colored faces }of its associated colored
cone. This bijection is constructed as follows:

Let  $Y\subset X$ be a  $G$-orbit. We denote by  ${\mathcal
  C}_Y(X)$  the cone generated by the image by $\rho$ of the
$B$-stable irreducible  divisors  which contain  $Y$\,. If ${\mathcal F}(Y)$
is the set of colours such that their closure contains  $Y$\,, then 
$({\mathcal C}_Y(X), {\mathcal F}(Y))$ is the colored face we are
looking for.

{\em iii) } If $P\subset G$ is the parabolic subgroup given by the
intersection of the stabilizers of the colours in $\mathcal D\setminus
\mathcal F(Y)$\,, then
\[
\dim{Y}=\operatorname{rg}(G/H) -\dim{\mathcal C_Y(X)}+\dim{G/P}\ .
\]

In particular, the unique closed orbit is a point if and only if the
dimension of $\mathcal C(X)$ is equal to the rank of $G/H$\,, and the
set of colours of the embedding is all ${\mathcal D}$\,. 

{\em iv) } A spherical variety $X$ is affine if and only if $X$ is
simple and there exists $\chi\in \Lambda_{G/H}$ such that

\begin{itemize}
\item[(a)] $\chi|_{{\mathcal V}(G/H)}\leq 0$

\item[(b)] $\chi|_{\mathcal C(X)} =0$

\item[(c)] $\chi|_{\rho ({\mathcal D(G/H)}\setminus {\mathcal F(X)})
    } > 0$\,. 
\end{itemize}
\qedennonce
\end{prop}

In  \cite{kn:vust}, Vust  established the combinatorial data
associated to symmetric spaces, in particular to the homogenous space 
$G\cong (G\times G)/\Delta(G)$\,, where  $G$ is a reductive connected
group, and $\operatorname{car}k=0$\,.
In order to carry out the classification of reductive monoids, we must
generalize his results to arbitrary characteristic. This is
essentially done in the literature, so we will state the results
without proving them. For a proof, we refer the reader to \cite{kn:tesis}.

Let us fix the notation:

We note $G_0$ the commutator of $G$\,. Thus, $G=ZG_0$\,, where
$Z=Z(G)$ is the connected center of $G$\,.  Moreover,  $Z_0=G_0\cap Z$
is finite; $G$  is then  isomorphic to $( G_0\times Z)/Z_0$ :
\[
k[G]\cong k[(G_0\times Z)/Z_0]\cong (k[G_0]\otimes_k
k[Z])^{Z_0}\,.
\]
 
We fix $T_0$\,, a  maximal torus of $G_0$\,, and 
$B_0\subset G_0$\,, a  Borel subgroup that contains  $T_0$\,.
Let $\alpha_1,\ldots,\alpha_l$ and
$\omega_1,\ldots ,\omega_l$ be the simple roots  and fundamental
weights associated to  $(B_0,T_0)$ respectively. 
We set $T=T_0Z$\,, a maximal torus of $G$\,. We note $W$ the
Weyl group associated to  $(G,T)$\,, and  $C=C(G)$ the
Weyl chamber of  $G$ associated to $B$\,.

If $B=B_0Z$\,, the weight lattice of   $G$ associated to 
$(B,T)$\,, satisfies:
\[
\Xi(T)\cong (\Xi(T_0)\times \Xi(Z))^{Z_0} = \{(\lambda,\mu)\in
\Xi(T_0)\times  \Xi(Z) \ :\ \mu|_{Z_0}=\lambda|_{Z_0}\}\,.
\]

\noindent It is a  sub-lattice of finite index of
$\Xi(G_0)\times \Xi(Z)$\,.

The  reflection associated to a simple root $\alpha$ is noted as usual
$s_\alpha$\,.

From the Bruhat decomposition, it is easy to prove the following 

\begin{prop}\label{Q(G)}{\rm (\cite{kn:tesis})}
Let $G$ be a reductive group. If we keep the preceding notations and
identify   $\operatorname{Hom}_\Z(\Xi(T),\Z)$ with $\Xi_*(T)$\,, then
\[
\begin{array}{ll}
{\mathcal Q}((G\times G)/\Delta(G)) = &{\mathcal Q}(G) = \Xi_*(T)\otimes
\Q\cong \vspace*{8pt} \\
 &(\Xi_*(T_0)\otimes \Q)\oplus  (\Xi_*(Z)\otimes\Q )=\mathcal
Q(G_0)\oplus \mathcal Q(Z)\ ,
\end{array}
\]

\noindent where the lattice considered is  $(\Xi_*(T_0)\times
\Xi_*(Z))^{Z_0} \subset\Xi_*(T_0)\times \Xi_*(Z)$\,. 

We identify  $\Xi_*(T)$ with $\Xi(T)$ by means of a $W$-invariant form
$\langle\cdot,\cdot\rangle$\,, in such a way that
$\langle \omega_i,\alpha_j^\vee \rangle=\delta_{ij}$\,. Under this
identification, $\mathcal Q(G)=\Xi(T)\otimes \Q$\,,  $\mathcal
Q(G)^*=\Xi(T)\otimes \Q$\,, and the valuation cone $\mathcal V(G)\subset
\mathcal Q(G)$ is  opposite to $C(G)$\,.

Moreover, the colours of  $G$ are the divisors 
$D_i=\overline{Bs_{\alpha_i}B^-}$\,, $i=1,\ldots l$\,, with associated
valuations $(\alpha_i^\vee,0)\in(\Xi(G_0)\times \Xi(Z))^{Z_0}$\,.
\qedennonce
\end{prop} 

It follows immediately that the behaviour 
under quotients by a finite central subgroup  of the combinatorial  data
associated to reductive monoids 
is ``good'':

\begin{prop}\label{quotient_par_groupe_central} {\rm
(\cite{kn:tesis},\cite{kn:vinberg}) }
Let  $G$ be a reductive  group and  $\Gamma \subset G$ a finite
central subgroup. Then
\ \\
\vspace*{-12pt}
\begin{itemize}
\item[i)] $\mathcal Q(G/\Gamma)=\mathcal Q(G)$\,, where we consider
the lattice $\Xi_*(T/\Gamma)=\Xi_*(T)^\Gamma\subset \Xi_*(T)$\,.

\item[ii)] If $\pi:G\rightarrow G/\Gamma$ is the canonical projection,
then $\pi$ induces a  bijection between $\mathcal D(G)$ and $\mathcal
D(G/\Gamma)$\,. 

\item[iii)] Under the identification i), $\rho_G(\mathcal D
(G))=\rho_{G/\Gamma}(\mathcal D(G/\Gamma))$\,, and  the valuation
cones $\mathcal V(G)$ and $\mathcal V(G/\Gamma)$ coincide. In
particular, $\pi_*$ is the identity. 

\item[iv)] If  $S$ is a  monoid with unit group $G$\,, then $S
/\Gamma$ is a reductive monoid, with unit group
$G/\Gamma$\,. Moreover, if we denote $\pi:S\rightarrow S/\Gamma$ the
quotient, then $\pi_*:\mathcal Q(G)\rightarrow \mathcal
Q(G')$  is the identity, and  under this identification we have ${\mathcal
C}(S/ \Gamma)={\mathcal C}(S)$\,.  
\end{itemize}
\qedennonce
\end{prop}

\begin{rem}
In the case of reductive embeddings, the application $\rho:\mathcal
D\rightarrow \mathcal Q(G)$ is injective. We identify then 
 $\mathcal D$ with its image by $\rho$\,.
\end{rem}

From this, we can recover the fact that a reductive monoid is a quasi-direct
product of a reductive group and a reductive monoid with zero
(\cite{kn:putcha-monoids},\cite{kn:tesis}): 

\begin{prop} \label{prop:quasi-direct_product}
Let $S$ be a reductive monoid with unit group $G$\,. Then there exist a
reductive group $G_1$ and a reductive monoid with zero $S_0$ such that
\[
S\cong (G_1\times S_0)/\Gamma\,,
\]

\noindent  where $\Gamma$ is a finite central
subgroup of $G$\,. In particular, if $G_2=G(S_0)$\,, then $G\cong (G_1\times
G_2)/\Gamma$\,.  
\qedennonce
\end{prop}

It follows that, as in the case of toric varieties, we can restrict
ourselves to the study of reductive monoids with zero:

\begin{prop}\label{passage_monoides_zero}
Let $S$ be a reductive  monoid with unit group $G$\,. Let $G_1,G_2, S_0$
and $\Gamma$ be as above. Then, under the identification $\mathcal
Q(G)=\mathcal Q(G_2)\oplus \mathcal Q (G_1)$\,, we have: 
\begin{equation}\label{equation_cone}
\left(\mathcal C(S),\mathcal F(S)\right)=\left(\{0\}\times \mathcal
  C(S_0), \{0\}\times \mathcal D(G_1) \right)\ .
\end{equation}
\end{prop} 

\dem
 By Propositions \ref{prop:quasi-direct_product} and
\ref{quotient_par_groupe_central} , we can suppose without loss of
generality that 
$G=G_1\times G_2$\,, and $S=G_1\times S_0$\,, $S_0$ a monoid with zero
and with unit group $G_2$\,. In particular,  $\Xi(G)=\Xi(G_1)\times
\Xi(G_2)$\,, and $\mathcal Q(G)=\mathcal Q(G_1)\oplus \mathcal Q(G_1)$\,.
If $B_i\subset G_i$ is a  Borel subgroup of $G_i$\,,
$i=1,2$\,, then $B=B_1\times B_2$ is a Borel subgroup of $G$\,.
If $\alpha_1,\cdots \alpha_r$ are the simple roots associated to
$B_1$\,, and $\beta_1,\ldots, \beta_s$ those associated to
$B_2$\,, then 
$(\alpha_1,0),\ldots,(\alpha_r,0),(0,\beta_1),\ldots,(0,\beta_s)$
are the simple roots associated to $B$\,.

Denote  $D_{\alpha_i} $ and $D_{\beta_j}$ the colours of $G_1$ and
$G_2$ respectively. It is immediate that the colours of $G$ are:
\[
\mathcal D(G)=\left\{ D_{\alpha_i}\times G_2\ ,\  G_1\times D_{\beta_j}\ 
  :\ i=1,\ldots r\ ,\ j=1,\ldots , s\right\}\ .
\]

As $k[G]=k[G_1]\otimes k[G_2]$\,, the valuations associated to
$D_{\alpha_i}\times G_2$  and  $G_1\times D_{\beta_j}$ are
$(\alpha^\vee_i,0)$ and 
$(0,\beta_j^\vee)$ respectively. In the same way, we see that
$\mathcal V(G)=-C(G)=-(C(G_1)\times C(G_2))$\,.  

A divisor  $D\subset S$ is $(G\times G)$-stable if and only if
$D=G_1\times D'$\,,  where $D'$ is  a divisor $(G_2\times
G_2)$-stable of  $S_0$\,. It follows that the valuation associated to
$D$ satisfies $\rho_G(\nu_D)=(0,\rho_{G_2}(\nu_{D'}))$\,. 

On the other hand, the kernel of  $S$ (i.e. the unique $(G\times
G)$-closed orbit) is equal to $G_1\times \{0\}$\,. If follows that the
set of colours $\mathcal F(S)$ is  $\{0\}\times \mathcal D(G_2)$ (we
recall that the set of colours of a monoid with zero $S_0$ is
$\mathcal D(G(S_0))$\,), and the proposition is proved.
\qqed

\begin{cor}\label{plongaffine}
Let  $S$ be a reductive monoid with unit group $G$\,. Let 
$({\mathcal   C}(S), {\mathcal F}(S))$  be the colored cone associated
to  $S$\,. Then ${\mathcal C}(S)+\Q^+{\mathcal D}$ is a strictly
convex cone.
\qedennonce
\end{cor}

The converse to this corollary is not true, as the following example
shows: 

\begin{example}
Consider $G=\operatorname{Sl}(2)\times k^*$\,. We identify
$\mathcal Q (G)$  to $\Q^2$\,, and  $\Xi_*(G)$ to the lattice
$\Z^2\subset \Q^2$\,. The affine embeddings of  $G$ are exactly those
where the closed orbit is a point (with colored cone of the form
 $(\mathcal C, (\alpha^\vee,0))$\,),
and the {\em elementary embeddings }associated to the half lines
$(\Q^+(0,1),\emptyset)$ and  $(\Q^+(0,-1),\emptyset)$\,.

Indeed, let  $(\mathcal C,\emptyset)$ be the colored cone associated
to an affine embedding without colors. Then there exists $\chi\in \Lambda$
such that 
$\chi(\alpha^\vee,0)>0$\,, $\chi(\nu)\leq 0$ for all $\nu\in \mathcal
V$\,. It follows that $\chi(\nu)=0$ if and only if 
 $\nu$ belongs to the  $y$-axis\,.  Let $\nu\in\mathcal V$ such that
the cone $\mathcal C'$ generated by $\nu$ and $(\alpha^\vee,0)$ is strictly
convex. Then the embedding $S$ associated to $(\mathcal C',(\alpha^\vee,0))$ is
affine (just take $\chi\equiv 0$\,). In particular, $S$ is a reductive
monoid (with zero: the set of colors of $S$ is all $\mathcal D$\,).

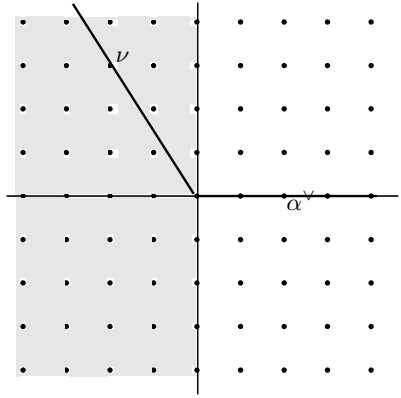
\begin{center}
\begin{picture}(0,0)%
\epsfig{file=figure1.pstex}%
\end{picture}%
\setlength{\unitlength}{0.00083300in}%
\begingroup\makeatletter\ifx\SetFigFont\undefined%
\gdef\SetFigFont#1#2#3#4#5{%
  \reset@font\fontsize{#1}{#2pt}%
  \fontfamily{#3}\fontseries{#4}\fontshape{#5}%
  \selectfont}%
\fi\endgroup%
\begin{picture}(2484,2489)(289,-1783)
\put(2051,-609){\makebox(0,0)[lb]{\smash{\SetFigFont{9}{14}{\rmdefault}{\mddefault}{\updefault}$\alpha^\vee$}}}
\put(985,314){\makebox(0,0)[lb]{\smash{\SetFigFont{9}{4}{\rmdefault}{\mddefault}{\updefault}$\nu$}}}
\end{picture}
\end{center}

On the other hand, if  we take $\tau=\mathcal C'\cap \mathcal V$\,,
then $\tau+\Q^+(\alpha^\vee,0)=\mathcal C'$ is strictly convex, but
the embedding associated to the colored cone $(\tau,\emptyset)$ is not
affine.
\end{example}

An immediate consequence of the preceding proposition is that the
colored cone of a reductive monoid $S$ is determined by $\mathcal V(S)$\,,
or even $\mathcal 
B(S)$\,. Moreover, we can find the decomposition of $S$ as a
quasi-direct product of a reductive monoid with zero and a
reductive group  (cf. Proposition \ref{prop:quasi-direct_product}) as follows:

If we keep  the notation of  Proposition \ref{prop:quasi-direct_product},
the space  $\mathcal Q(G_2)$ is $\langle \mathcal C(S)\rangle$\,, the
sub-space of $\mathcal Q(G)$ generated by $\mathcal C(S)$\,. The
lattice considered contains $\Xi(G)\cap \langle
\mathcal C(S)\rangle$ as a subgroup of finite index.

A simply connected covering of the commutator of the normal
subgroup $G_1$ (resp. $G_2$\,) is the semi-simple simply connected group
whose  root system is given by the  Dynkin diagram
of $G$ restricted to $\mathcal F(S)$ (resp. $\mathcal 
D(G)\setminus \mathcal F(S)$\,).

In terms of the combinatorics of the colored cones of  $\mathcal
Q(G)$\,, this remark gives the following result:

\begin{prop}\label{construction_monoides}
Let $G$ be a reductive group. Given a polyhedral cone $\tau$ contained
in ${\mathcal V}(G)=-C$\,, such that  $\tau+\Q^+{\mathcal D}$ is
strictly convex, there exists an unique subset of colours  ${\mathcal 
F}\subset {\mathcal D}$ such that  $(\tau +\Q^+{\mathcal F}, {\mathcal
  F})$ is the colored cone associated to an affine embedding of $G$\,.
\end{prop} 

\dem
Let $\tau\subset {\mathcal V}(G)$ be a  cone such that 
$\tau+\Q^+{\mathcal D}$ is strictly convex. 
If  $A$ is a subset of $\mathcal Q(G)$\,, denote 
 $\langle A\rangle $ the sub-space of $\mathcal Q(G)$ generated by  $A$\,.
Consider the  projection $p: {\mathcal Q}(G)\rightarrow {\mathcal
  Q}(G_0)$\,. 

By the preceding remark,  $\mathcal F=\mathcal D\cap
\langle p(\tau)\rangle$ is the natural candidate for  $(\tau+
\Q^+\mathcal F,\mathcal F)$ to be a colored cone correponding to an
affine embedding of $G$\,.  Let  $\{v_1,\ldots,
v_r\}\subset \mathcal 
V(G)$ be generators of the cone $\tau$\,. Then  $\tau+\Q^+{\mathcal F}$ is
generated by  $\{v_1,\ldots, v_r\}\cup {\mathcal F}\subset {\mathcal
  V}(G)\cup {\mathcal D}$\,. The cone  $\tau+\Q^+{\mathcal D}$ is
strictly convex, and thus  the same holds for $\tau+\Q^+{\mathcal F}$ . 

Let us show that  $(\operatorname{int}(\tau+\Q^+ \mathcal
F))\cap\mathcal V(G)$ is non-void. Then $(\tau+\Q^+\mathcal F,\mathcal
F)$ will be a colored cone.

Let $v=\sum_{n_j>0} n_jv_j\in
\operatorname{int}{\tau}$\,. We have: 
\[
p(v)\in 
p\left(\operatorname{int}\tau\right)\subset p({\mathcal V}(G))=-C(G_0) \ .
\]

If  $q>0$ is sufficiently small, as the  set  ${\mathcal F}$ is
contained in $ \langle p(\tau)\rangle$\,, we have 
\[
p\left(v+\sum_{(\alpha^\vee,0)\in {\mathcal F}}q(\alpha^\vee,0)\right) =
p(v)+\sum_{(\alpha^\vee,0)\in {\mathcal F}}q(\alpha^\vee,0)\in
\operatorname{int}(p(\tau))\,.
\]

It follows that
\[
v+q\sum_{(\alpha^\vee,0)\in {\mathcal F}}(\alpha^\vee,0) \in  
(\operatorname{int}(\tau+\Q^+{\mathcal F}))\cap p^{-1}(-C(G_0))=
(\operatorname{int}(\tau+\Q^+{\mathcal F})) \cap {\mathcal V}(G)\,, 
\]

\noindent and thus $(\tau + \Q^+{\mathcal F}, {\mathcal F})$ 
is a colored cone. 

In order to prove that $(\tau + \Q^+{\mathcal F}, {\mathcal F})$
corresponds to an affine embedding, we consider a  basis ${\mathcal
F}\cup \{w_1,\cdots w_s\}$ of   $\langle \tau \rangle +\langle
{\mathcal    F}\rangle$\,, where  $w_i$ belongs to  $\Xi(G)$ for all
$i=1,\ldots s$\,. 
By construction, $\left( \langle \tau \rangle +\langle {\mathcal
F}\rangle\right) \cap \langle \mathcal D\setminus {\mathcal F}\rangle
= \{0\}$\,. We can then add the elements of ${\mathcal
D}\setminus \mathcal F$ to $\mathcal F\cup \{w_1,\ldots,w_s\}$ and
obtain in this way a linearly independent set. We complete this set in
a basis $\mathcal B$ of  ${\mathcal Q}(G)$\,:
\[
{\mathcal B}={\mathcal D}\cup \{w_1,\ldots, w_s\}\cup
\{w_{s+1},\ldots, w_{h}\} \quad \quad w_i\in \Xi(G)\ ,\ i=1,\ldots, h\  .
\]

Define  $\chi\in \Xi_+(G)$ over  ${\mathcal B}$ as follows:
\[
\chi(b)=\left\{
\begin{array}{ll}
1\quad\quad & {\rm if}\ b\in  \mathcal D \setminus {\mathcal F}\\
&\\
0& {\rm otherwise}
\end{array}
\right.
\]

\noindent and extend $\chi$ by linearity. It is easy to see that
$\chi$ satisfies the  conditions a)--c) of 
Proposition \ref{bijection_orbites_faces} (iv), so the embedding associated
to  $(\tau +\Q^+{\mathcal F},{\mathcal F})$ is affine.

It remains to prove the unicity of the chosen set of colours.
Let   $\mathcal F'\subset \mathcal D$ be a subset of colours such that
 $(\tau+\Q^+\mathcal F', \mathcal F')$ is a colored cone corresponding
to an affine embedding of $G$\,. Let us show that  $\mathcal F'=\mathcal F$\,.

Let  $v\in (\operatorname{int}(\tau +\Q^+\mathcal F')) \cap \mathcal V(G) =
(\operatorname{int}\mathcal C(S'))\cap \mathcal V(G) $\,. The vector
$v$ can be written as
\[
v = v_0 + \sum_{(\alpha,0)\in \mathcal F'} q_\alpha(\alpha^\vee,0)\quad
  ,\quad v_0\in \tau  \,,\ q_\alpha>0 \ .
\]

Let  $\chi\in \Lambda$ be as above.  As $\chi(w)\leq 0$ for all  $w\in
\mathcal V(G)=-C(G)$\,, and $\chi>0$ over $\mathcal D\setminus 
\mathcal F'$\,,  it follows that  $0\geq \chi(v)=\sum_{(\alpha,0)\in
\mathcal F'\setminus   \mathcal F} q_\alpha\chi(\alpha^\vee,0)$\,.
The only possibility is that $\chi(v)=0$\,, and 
$\mathcal F'\subset \mathcal F$\,. By symmetry, we get $\mathcal
F\subset \mathcal F'$\,.
\qqed

\begin{cor}\label{corolario}
Let  $G$ be a  reductive group. Let  $S$ and  $S'$ be two monoids with
unit group $G$\,, such that $\Q^+\mathcal B(S)=\Q^+\mathcal
B(S')$\,. Then  $S$ and  $S'$ are isomorphic. In particular, if $
\mathcal V(S)=\mathcal V(S')$\,, then $S$ and $S'$ are isomorphic.
\qedennonce
\end{cor}

\begin{rem}
Let us give a geometric interpretation of the corollary above. 
If $S$ is a reductive monoid with unit group $G$ and 
$(\mathcal C(F),\mathcal F(S))$ its colored cone, there exists a
proper $G$-morphism $\widetilde{S}\rightarrow S$\,, where
$\widetilde{S}$ is the embedding associated to $(\mathcal
V(S),\emptyset)$ (\cite{kn:brion-geovarsph}). Moreover, the embedding
$\widetilde{S}$ is minimal for this property; it is called the {\em
decoloration }of $S$\,. Corollary 
\ref{corolario} says then that two monoids are isomorphic if and only
if their decolorations are.
\end{rem}

In conclusion, we have proved the following classification:

\begin{teo}
Let $G$ be a reductive group. The isomorphism classes of algebraic
monoids with unit group $G$ are in bijection with the strictly
convex polyhedral cones of $\mathcal Q(G)$ generated by all the colours and a
finite set of elements of $\mathcal V(G)$\,.
\end{teo}

\section{Some consequences}

We end up with some applications of the classification obtained above.

First, we describe the algebra of regular functions on a reductive
monoid,  generalising the result obtained by Vinberg in
characteristic zero (\cite{kn:vinberg}).

\begin{nota}
Let $G$ be a reductive group. We keep the notations of the last
 section. If  $M$ is a $G$-module, we note $M^{(B)}$ the set of
 eigenvectors of $M$ for the action of $B$\,.
\end{nota}

\begin{prop}\label{qwert}
Let $S$ be a monoid with unit group $G$\,. Then,
\[
k\left[S\setminus \overline{\cup_{D\in \mathcal F(S)}D}\right]^{(B\times
  B^-)}=\left\{ f\in 
  k(G)^{(B\times B^-)}\ :\ \chi_f\in {\mathcal 
C}(S)^\vee \right\}\ , \quad and
\]
\[
k[S]^{(B\times B^-)}=\left\{ f\in 
  k(G)^{(B\times B^-)}\ :\ \chi_f\in {\mathcal 
C}(S)^\vee \cap \Xi_+(G)\right\}\ .
\]

Moreover, the set  $\mathcal L_S=\mathcal C(S)^\vee\cap \Xi(T)$ of 
weights of $k[S]^{(B\times B^-)}$ generates $\Xi(T)$\,, the set of
weights of $k(G)^{(B\times B^-)}=k(S)^{(B\times B^-)}$\,.

In particular, if $\operatorname{car}k=0$\,, the decomposition of
$k[S]$ into simple $(G\times G)$-modules is the following: 
\[
k[S]=\bigoplus_{ \chi\in {\mathcal C}(S)^\vee\cap \Xi_+(G)}
V_\chi\otimes V_\chi^*\ .
\]
\end{prop}

\dem
The first equality follows from \cite[thm. 3.5]{kn:knop}. In order to
prove the second one, observe that a function  $f\in k\left[S\setminus
\overline{\cup_{D\in \mathcal  F(S)}D}\right]^{(B\times B^-)}$
extends to a regular function on $S$ if and only if
$\nu_D(f)\geq 0$ for all $D\in \mathcal F(S)$\,. By proposition
\ref{Q(G)}, this means that $\langle 
\chi_f,\alpha_i^\vee\rangle=\nu_{D_i}(f)\geq 0$\,, i.e.
$\chi_f$ is a dominant weight.

Let $f=\psi/\varphi\in k(S)^{(B\times B^-)}$\,, where $\varphi,\psi\in
k[S]$\,. The  $(G\times G)$-module $V\subset k[S]$ generated by
$\varphi$ is of finite dimension. Let $\{g_1\cdot \varphi,\ldots
,g_r\cdot\varphi\}$\,, $g_i\in G\times G$\,, be a system of generators
of $V$\,. Then, there exists a $(B\times B^-)$-eigenvector
$h=\sum_{i=1}^{r} t_ig_i\cdot \varphi\in V$\,, $t_i\in k$\,. Consider
the $(B\times B^-)$-eigenvector $hf=(\sum_{i=0}^r t_i g_i\cdot
\varphi)\psi/\varphi\in k[S]$\,; its weight is
$\chi_{hf}=\chi_h+\chi_f$\,. We have  $f= hf/h$\,, where $hf,h\in
k[S]^{(B\times B^-)}$\,, and $\chi_f=\chi_{fh}-\chi_h$\,. It follows
that $\mathcal L_S$ generates $\Xi(T)$\,.

Finally, the last equality is a direct consequence of the complete
reducibility of the rational representations of $G\times G$ in
characteristic $0$\,.
\qqed

By Proposition \ref{prop:ideals-idempotents}, the ideals of a
reductive monoid  $S$ are
exactly the closed 
$(G\times G)$-stables subvarieties of $S$\,. The dimension of an
irreducible component of an ideal is then given by the dimension of the
open orbit in this component. By  Proposition
\ref{bijection_orbites_faces} (iii), to calculate those dimensions it
suffices to 
calculate the stabilizers of the colours. But the Bruhat decomposition
gives us this information immediately:

\begin{prop}[{\cite[\S 29]{kn:humpreys}}]\label{lem:stabilisateurs-couleurs}
Let $G$ be a reductive group. If we keep the preceding notations, the
stabilizer in  $G\times G$ of the colour
$D_i=\overline{Bs_{\alpha_i}B^-}$\,, $i=1,\ldots ,l$\,, is generated by the
parabolic subgroups $(P_{\alpha_j}\times P_{-\alpha_j})$\,,  $j\neq i$\,,
where $P_{\alpha}$ denotes the minimal parabolic associated to the
simple root  $\alpha$\,. 
\qedennonce
\end{prop}

\begin{cor}\label{dimension_orbites}
Let  $S$ be a reductive monoid with unit group $G$\,. If $Y\subset S$ is a
$(G\times G)$-orbit, whose combinatorial data is $({\mathcal C}_Y(S),{\mathcal
  F}(Y))$\,, its dimension is given by:
\[
\dim Y= l+\dim Z(G) -\dim {\mathcal C}_Y(S) + \dim G/P_{\mathcal
  D\setminus \mathcal F(Y)}
\]

\noindent where $P_{\mathcal D\setminus \mathcal F(Y)}$ is the 
 parabolic subgroup generated by the parabolics $P_{\alpha}$\,, $\alpha\in
\mathcal D\setminus \mathcal F(Y)$\,. 
\end{cor}

\dem

Apply Proposition \ref{bijection_orbites_faces} (iii)  and then Lemma
\ref{lem:stabilisateurs-couleurs} in order to calculate the intersection of
the parabolic subgroups which stabilize the colours belonging to
$\mathcal D\setminus\mathcal F(Y)$\,.
\qqed

\vspace*{1cm}
\noindent {\sc Alvaro Rittatore}\\
Institut Fourier,\\
Grenoble I.\\
B.P. 74.\\
38402 St. Martin d'H\`eres CEDEX
France\\
e-mail: {\tt alvaro@fourier.ujf-grenoble.fr}

\begin{thebibliography}{99}

\bibitem{kn:brion-geovarsph}  M. Brion, {\em Sur la G\'eom\'etrie des
vari\'et\'es sph\'eriques. }Comment. Math. Helv. 66 (1991), 237--262.
\bibitem{kn:hartshorne-amplesubvarieties} R. Hartshorne. {\em Ample
subvarieties of Algebraic Varieties. } Lect. Notes in Math
157. Springer Verlag, New York, 1970.
\bibitem{kn:humpreys} J. Humpreys, {\em Linear Algebraic Groups.} GTM
21, Springer Verlag, New York, 1970.
\bibitem{kn:kempf} G. Kempf, {\em Algebraic Varieties. } London
Math. Soc. Lect. Notes 172, Cambridge Univ. Press, 1993.
\bibitem{kn:knop} F. Knop, {\em The Luna-Vust Theory of Spherical
Embeddings. }In 
S. Ramanan {\em et al}, editors, Procedings of the Hydebarad
Conference on Algebraic Groups, pages 225--249. National Board for
Higher Mathematics, Manoj, 1991.  
\bibitem{kn:mumford} D. Mumford, J. Fogarty and F. Kirwan, {\em
Geometric Invariant Theory. } Modern Surveys in Math. 34, Springer
Verlag, New York, 3rd enlarged ed., 1994.
\bibitem{kn:oda} T. Oda, {\em Torus embeddings and applications. }
Tata Press, Bombay, 1978.
\bibitem{kn:putcha-monoids} M.S. Putcha,  {\em Linear algebraic
monoids. }London Math. Soc. Lect. Notes 133, Cambridge Univ. Press, 1988.
Notes 133,  
\bibitem{kn:renner-thesis} L. Renner, {\em Algebraic monoids. }
PhD. Thesis, Univ. of British Columbia, 1982.
\bibitem{kn:tesis} A. Rittatore, {\em Mono\"\i des alg\'ebriques et
vari\'et\'es sph\'eriques. }PhD. Thesis, Institut Fourier, Grenoble, France,
1997. {\tt http://www-fourier.ujf-grenoble.fr} 
\bibitem{kn:vinberg} E.B. Vinberg, {\em On reductive Algebraic
Semigroups. }Amer. Math. Soc. Transl., Serie 2, 169 (1994), 145 --
182. Lie Groups and Lie Algebras. E.B. Dynkin seminar.
\bibitem{kn:vust} T. Vust, {\em Plongements d'espaces sym\'etriques
alg\'ebriques: une  classification. }Ann. Scoula Norm. Sup. Pisa,
XVII, 2 (1990), 165--194. 
\bibitem{kn:waterhouse} W. Waterhouse, {\em The unit group of affine
algebraic monoids. }Proc. Amer. Math. Soc. 85 (1982) 506 -- 508.

\end{thebibliography}
\end{document}